\documentclass[11pt]{article}
\usepackage[dvips]{color}
\usepackage{graphics}
\usepackage{amsmath}
\usepackage{amsfonts}
\usepackage{latexsym}
\usepackage{amssymb}
\newcommand{\g}{{\mathfrak g}}
\newcommand{\h}{{\mathfrak h}}
\newcommand{\n}{{\mathfrak n}}

\newcommand{\C}{{\bf C}}
\newcommand{\R}{{\bf R}}

\newcommand{\Z}{{\bf Z}}

\usepackage{multirow}

\newsavebox{\savepar}

\begin{document}

\title
{Complex and K\"ahler structures on compact homogeneous
manifolds - their existence, classification and moduli problem\\}

\author{Keizo Hasegawa\\
}
\date{
}
\maketitle

\bigskip
%\vspace{0.25in}
\noindent{\large\em \bfseries 1.  Introduction}
\bigskip

\noindent In the field of complex geometry one of the primary problems is
whether given real manifolds admit certain complex geometrical
structure such as complex structures, K\"ahler structures, or Stein
structures.  For the case of complex structures, for instance, we have a long standing
problem of whether $S^6$ admits a complex structure.
We can consider this problem in a more general setting: 
we extend the problem to the case of homogeneous
manifolds of compact semi-simple Lie groups, including $S^6$ and ${\bf C}P^n$;
and propose a closely related conjecture:
\smallskip

\noindent {\bf Conjecture 1}. 
{\em A compact homogeneous manifold of compact semi-simple Lie group admits
only homogeneous complex structures.}
\smallskip

We see, according to Wang's classification of compact simply connected homogeneous
complex manifolds [\ref{W}], that $S^6$ admits no homogeneous complex structures; and thus
the above conjecture implies that
$S^6$ admits no complex structures. We also know that the complex structure of ${\bf C}P^n$
is {\em rigid}, that is,  admits no non-trivial small deformations; and ${\bf C}P^2$ admits no
complex structures but the original one (which is homogeneous complex). More generally,
we have a class of flag manifolds including ${\bf C}P^n$, which are simply connected,
homogeneous K\"ahlerian; and
according to Borel [\ref{B}], they are the only homogeneous manifolds of compact semi-simple
Lie groups which admit K\"ahler structures. On the other hand, according to Samelson [\ref{S}],
compact semi-simple Lie groups of even dimension (being considered as homogeneous
manifolds) admit homogeneous complex structures but no K\"ahler structures (since $b_2=0$).
It seems that there are not yet known any counter-examples to the above conjecture. 
\medskip

For the case of compact homogeneous manifolds of dimension 4, we have a complete
classification of those which admit complex structures (see Section 2 and [\ref{H2}]);
in particular we showed that
any  complex structure on a $4$-dimensional compact solvmanifold
$\Gamma \backslash G$ (up to finite covering) is {\em left-invariant}, that is, induced from
a left-invariant complex structure on a simply connected solvable
Lie group $G$ with lattice $\Gamma$. Furthermore, we recently showed (see Section 5 and
[\ref{H3}]),
based on Nakamura's results [\ref{N}], that there exists a compact 
solvmanifold of dimension 6---actually a compact complex solvmanifold of dimension 3---which
admit a continuous family of non-left-invariant complex structures
as small deformations of  the original complex structure.
This result is important since it implies
that there are \char96\char96 abundant" complex structures on compact solvmanifolds of higher
dimension, while as mentioned above, it seems that compact homogeneous manifolds of
compact semi-simple Lie groups admit only \char96\char96 restricted" complex structures.
Concerning left-invariant complex
structures on simply connected solvable Lie groups, we have a related conjecture:
\smallskip

\noindent {\bf Conjecture 2}. 
{\em Any left-invariant complex structure on a $2n$-dimensional simply connected unimodular
solvable (nilpotent) Lie group is Stein (biholomorphic to ${\bf C}^n$ respectively)}.
\smallskip 

It should be noted that an $n$-dimensional simply connected complex solvable Lie group is
biholomorphic to ${\bf C}^n$ (cf. [\ref{N}]),  and the conjecture holds for $n=2$ [\ref{OR}]. 
\medskip

We have more decisive results for the classification problem of K\"ahler structures on
compact homogeneous manifolds. For the case of compact solvmanifolds, 
we have the following result:
\smallskip

\noindent {\bfseries Theorem} ([\ref{H2}, \ref{H3}]).
{\em A compact  solvmanifold admits a K\"ahler structure if and only if 
it is a finite quotient of a complex torus which has a structure of a complex torus bundle
over a complex torus. }
\smallskip

We can express a class of compact K\"ahlerian solvmanifolds in the theorem explicitly 
as those of the form $\Gamma \backslash G$, where $G$ is a simply connected 2-step solvable
Lie group with lattice $\Gamma$---they are exactly hyperelliptic surfaces for dimension 4
(see Section 2). For the case of reductive Lie groups,  the classification of compact
homogeneous K\"ahler manifolds suggests the following conjecture to us:
\smallskip

\noindent {\bf Conjecture 3}. 
{\em A compact homogeneous manifold of reductive Lie group admits a K\"ahler structure
if and only if it is the product of a complex torus and a flag manifold.}
\medskip

As mentioned before, a homogeneous manifold of compact semi-simple Lie group admits
a K\"ahler structure if and only if it is a flag manifold. On the other hand, we know that
$S^1 \times \Gamma \backslash \widetilde{\rm SL}_2({\bf R})$,
where $\widetilde{\rm SL}_2({\bf R})$ is the universal covering of
${\rm SL}_2({\bf R})$ (which is a simply connected non-compact semi-simple Lie group),
and $\Gamma$ is a lattice of $\widetilde{\rm SL}_2({\bf R})$, admits
a non-K\"ahler complex structure---defining an elliptic surface [\ref{W}].
\medskip

A pseudo-K\"ahler structure is a pseudo-Hermitian structure with its associated fundamental
form $\omega$ being closed (see Section 2).
It is known that a compact homogeneous pseudo-K\"ahler manifold is biholomorphic to the
product of  complex torus and a flag manifold [\ref{DG}]. There exists
a compact homogeneous complex pseudo-K\"ahler solvmanifold which is not homogeneous
pseudo-K\"ahlerian, that is, the pseudo-K\"ahler form $\omega$ may not be invariant by the group
action [\ref{Y1}].
The classification of compact homogeneous complex pseudo-K\"ahler solvmanifolds is not yet
known; at this moment, we have a complete classification for complex dimension 3
(see Section 4 and [\ref{H3}, \ref{H4}]), and a structure theorem that
it is  a holomorphic fiber bundle over a complex torus with fiber a complex torus [\ref{Y2}], 
which is, unless trivial, not a principal bundle. Recently Guan [\ref{Gu}] has shown a
fundamental theorem on cohomology groups of compact solvmanifolds,
which could be applied to our problem.
\bigskip

%%%%%%%%%%%%%%%%%%%%%%%%%%%%%%%%%%%%%%%%%%
\noindent{\large\em \bfseries 2.  Complex and K\"ahler structures on
compact homogeneous manifolds}
\bigskip

%\fbox{\begin{minipage}{18cm}
\noindent{\em A homogeneous manifold} $M$ is a differentiable manifold
on which a real Lie group $G$ acts transitively.
$M$ is a {\em homogeneous complex manifold}, if $M$ is a complex
manifold and the group action is holomorphic.
\smallskip

We will first make some important remarks:
%\end{minipage}
%}
\smallskip

\begin{list}{}{\topsep=5pt \leftmargin=15pt \itemindent=5pt \parsep=0pt \itemsep=5pt}
\item[ (1)] In the case where $M$ is a compact homogeneous complex manifold,
we can assume that $G$ is a complex Lie group [\ref{W}].
\item[ (2)] A Lie group $G$, as a homogeneous manifold, 
admits a homogeneous complex structure $J$
if and only if $J$ is a left-invariant complex structure on $G$.
\item[ (3)] A complex structure $J$ on a Lie group $G$ is both left and right-invariant
if and only if $G$ is a complex Lie group (w.r.t. $J$).
\end{list}
\bigskip

The complete classification of 2-dimensional compact homogeneous complex manifolds
are known:
\medskip

%%%%%%%%%%%Theorem 1
\noindent{\bfseries Theorem 1} (Tits [\ref{T}]).
%\fbox{\begin{minipage}{18cm}
{\em 2-dimensional compact homogeneous complex manifolds
are biholomorphic to one of the following:
\smallskip

$$T^2, \C P^2, \C P^1 \times \C P^1,
T^1 \times \C P^1, {\rm Homogeneous\; Hopf\; surface},$$
where $T^k$ denotes a $k$-dimensional complex torus.}
%\end{minipage}
%}
\medskip

For instance, we have 

\begin{list}{}{\topsep=5pt \leftmargin=15pt \itemindent=5pt \parsep=0pt \itemsep=3pt}
\item[ (1)]
$\C P^1 = B \backslash G = (B \cap H) \backslash H$,
where $G = {\rm SL}_2 (\C)$, $ H = {\rm SU}_2 (\C)$ and $B$ is a Borel subgroup of $G$:
$$B = \big\{ \left(
\begin{array}[c]{cc}
\alpha & \beta\\
0 & \delta\\
\end{array}
\right)
\in {\rm SL}_2 (\C)\big\}.
$$

\item[ (2)]
A homogeneous Hopf surface $S$ is by definition $W/ \Gamma_\gamma$,
where $W = (\C^2 - \{O\})$ and $ \Gamma_\gamma$ is a group of
automorphisms on $W$ generated by the multiplication by $\gamma\, (\not= 1)$.
$S$ is diffeomorphic to
${\rm U}_2 (\C) = {\rm SU}_2 (\C) \rtimes S^1 \cong S^3 \times S^1,$
and $S = B_\gamma \backslash G$, where $B_\gamma$ is the subgroup of $B$
with $\alpha = \gamma^k, \delta = \gamma^{-k}\, (k \in \Z)$.
Note that $S$ has a structure of a holomorphic $T^1$-bundle over $\C P^1$.
\end{list}
\bigskip

%\fbox{\begin{minipage}{18cm}
$M$ is a {\em homogeneous complex K\"ahler manifold},
if $M$ is a homogeneous complex manifold which admits
a K\"ahler structure.
$M$ is a {\em homogeneous K\"ahler manifold}, if $M$ is a homogeneous complex
manifold which admits  a K\"ahler structure invariant by the group action.
The following theorem (Theorem 2) is well known, which was first proved
by Matsushima for homogeneous K\"ahler cases, and later by Borel-Remmert
for homogeneous complex K\"ahler cases.
%\end{minipage}
%}

\bigskip
%%%%%%%%%%Theorem 2
\noindent{\bfseries Theorem 2} (Matsushima [\ref{M}], Borel-Remmert [\ref{BR}]).
%\fbox{\begin{minipage}{18cm}
{\em A compact homogeneous complex K\"ahler
manifold is biholomorphic to the product of a complex torus and a homogeneous rational
manifold (which is a compact simply connected algebraic manifold).}
%\end{minipage}
%}
\medskip

Let $M = \Gamma \backslash G$ be a compact homogeneous complex manifold,
where $G$ is a simply connected complex Lie group with discrete subgroup $\Gamma$.
Then, Theorem 2 implies that
$M$ admits a K\"ahler structure
if and only if $M$ is a complex torus.
In particular, the only compact homogeneous complex K\"ahler
solvmanifold is a complex torus.
\bigskip

%\fbox{\begin{minipage}{18cm}
%Let $M$ be a symplectic manifold with symplectic form $\omega$. If $M$
%admits a complex structure $J$ such that $\omega(JX, JY) = \omega(X, Y)$
%for any vector fields $X, Y$ on $M$, then we call $(\omega, J)$ a {\em pseudo-K\"ahler}
%structure on $M$. 

Let $M$ be a symplectic manifold with symplectic form $\omega$. If $M$
admits a complex structure $J$ such that $\omega(JX, JY) = \omega(X, Y)$
for any vector fields $X, Y$ on $M$, we call $(\omega, J)$ a {\em pseudo-K\"ahler}
structure on $M$. For a pseudo-K\"ahler structure $(\omega, J)$, we have a
pseudo-Hermitian structure $(g, J)$ defined by
$g(X, Y) = \omega(X, JY)$.
In other words, a pseudo-K\"ahler (K\"ahler) structure is 
a pseudo-Hermitian (Hermitian) structure $(g, J)$ with its closed fundamental form $\omega$,
where $\omega$ is defined by $\omega(X,Y) = g(JX, Y)$ for any vector fields $X, Y$.
\medskip

The following theorem (Theorem 3) asserts that Theorem 2 holds also for
homogeneous pseudo-K\"ahler cases.
However,  as we will see in Section 4,
it does not hold for compact homogeneous complex pseudo-K\"ahler 
manifolds.
%\end{minipage}
%}
\medskip

%%%%%%%%%Theorem 3
\noindent{\bfseries Theorem 3} (Dorfmeister and Guan [\ref{DG}]).
%\fbox{\begin{minipage}{18cm}
{\em A compact homogeneous pseudo-K\"ahler
manifold is biholomorphic to the product of a complex torus and a homogeneous rational
manifold.}
%\end{minipage}
%}
\bigskip

%\fbox{\begin{minipage}{18cm}
A compact solvmanifold $M$ can be written as, up to finite covering,
$$M = \Gamma \backslash G,$$
where $G$ is a simply connected real solvable Lie group and $\Gamma$ is a
lattice of $G$.
A complex structure $J$ on $M$ is {\em left-invariant complex structure}, if it
is deduced from a left-invariant complex structure on $G$.
%\end{minipage}
%}
\bigskip

A left-invariant complex structure
$J$ on $G$ can be considered as a linear automorphism of $\g$,
that is $J \in {\rm GL} (\g, \R)$,
such that $J^2 = - I$; and the integrability condition is satisfied:
$$N_J(X,Y)=[JX,JY]-J[JX,Y]-J[X,JY]-[X,Y]
$$
vanishes for $X, Y \in \g$.
\medskip

We have a complete list of complex structures
on  four-dimensional compact solvmanifolds, all of which are
left-invariant
[\ref{H2}]:
\medskip

%%%%%%%%%%%Theorem 4
\noindent{\bfseries Theorem 4} ([\ref{H2}]).
%\fbox{\begin{minipage}{18cm}
{\em A complex surface is diffeomorphic to a
four-dimensional solvmanifold if and only if it is one of the following surfaces:
Complex torus, Hyperelliptic surface, Inoue Surface of type $S^0$,
Primary Kodaira surface, Secondary Kodaira surface, Inoue Surface of type
$S^{\pm}$. 
Furthermore, every complex structure on each of these complex surfaces
(considered as solvmanifolds) is left-invariant.
} 
%\end{minipage}
%}
\bigskip

We can express each of these complex structures as a linear automorphism $J$ of $\g$.
In the following list,  for each surface the Lie algebra $\mathfrak g$ of $G$ has
a basis $\{X_1,X_2,X_3,X_4\}$ with only nonzero brackets specified.

Except for (6), the complex structure $J$ is defined by
$$JX_1=X_2, JX_2=-X_1, JX_3=X_4, JX_4=-X_3.$$

%list
\begin{list}{}{\topsep=0pt \leftmargin=5pt \itemindent=15pt \parsep=0pt \itemsep=10pt}
\item[ (1)] Complex Tori \par
$[X_i, X_j]=0 \;( 1 \le i < j \le 4);$
\item[ (2)] Hyperelliptic Surfaces \par
$[X_4, X_1]= -X_2, [X_4, X_2]= X_1;$
\item[ (3)] Inoue Surfaces of Type $S^0$ \par
$[X_4, X_1]= a X_1 - b X_2, [X_4, X_2]= b X_1 + a X_2,
[X_4, X_3]= -2a X_3,$
where $a, b \;(\not=0) \in \R$;
\item[ (4)] Primary Kodaira Surfaces \par
$[X_1, X_2]= -X_3;$
\item[ (5)] Secondary Kodaira Surfaces \par
$[X_1, X_2]= -X_3, [X_4, X_1]= -X_2, [X_4, X_2]= X_1;$
\item[ (6)] Inoue Surfaces of Type $S^+$ and $S^-$ \par
$[X_2, X_3]= -X_1, [X_4, X_2]= X_2, [X_4, X_3]= -X_3,$
and,
$$J X_1=X_2, J X_2=-X_1, J X_3= X_4-q X_2, J X_4= -X_3-q X_1.$$
\end{list}

\bigskip

%%%%%%%%Example 1
\noindent{\bfseries Example 1} (Hyperelliptic Surfaces).
Let $G =: (\C \times \R) \rtimes \R$, where the
action $\phi: \R \rightarrow {\rm Aut}(\C \times \R)$ is defined by 
$$\phi(t)((z,s)) = (e^{\sqrt{-1} \eta t} z, s),$$
where $\eta = \pi, \frac{2}{3} \pi, \frac{1}{2} \pi$ or $\frac{1}{3} \pi$.
\medskip

Since the action on the second factor $\R$ is trivial, the multiplication of $G$ is
defined on $\C^2$ as follows:
$$(w_1,w_2) \cdot (z_1,z_2)= (w_1+e^{\sqrt{-1} \eta t} z_1, w_2+z_2),$$
where $t= {\rm Re}\,w_2$.  
\medskip

We can see that there exist seven
isomorphism classes of lattices $\Gamma$ of $G$, which correspond to seven
classes of hyperelliptic surfaces.
\bigskip

%%%%%%%%%Example 2
\noindent{\bfseries Example 2} (Primary Kodaira Surfaces).
Let $G = N \times \R$ be the nilpotent Lie group, where 
$$N = \Bigg\{ \left(
\begin{array}[c]{ccc}
1 & x & s\\
0 & 1 & y\\
0 & 0 & 1
\end{array}
\right) \rule[-8mm]{0.25mm}{18mm}\;  x, y, s \in \R \Bigg\},
$$
which has the lattice $\Gamma_n$ with ${\displaystyle s = \frac{z}{n}}, x, y, z \in \Z$.

Taking the coordinate change $\Phi$ from $N \times \R$ to $\R^4$:
$$\Phi: ((x,y,s), t) \longrightarrow
(x,y, 2 s-x y, 2 t+ \frac{1}{2}(x^2+y^2)),$$
and regarding $\R^4$ as $\C^2$, the group operation on $G$ can be expressed as
$$(w_1,w_2) \cdot (z_1,z_2)=
(w_1 + z_1, w_2 - \sqrt{-1} \bar{w_1} z_1 + z_2).$$
\bigskip

%\fbox{\begin{minipage}{18cm}
Let $M$ be a solvmanifold of the form $ \Gamma \backslash G$.
$M$ is {\em of completely solvable type},
if the adjoint representation of $\mathfrak g$ has
only real eigenvalues.
$M$ is {\em of rigid type},
if the adjoint representation of $\mathfrak g$ has only
pure imaginary (including $0$) eigenvalues. We note that
%\end{minipage}
%}
\medskip

\begin{list}{}{\topsep=5pt \leftmargin=15pt \itemindent=5pt \parsep=0pt \itemsep=3pt}
\item[ (1)]
It is clear that $M$ is both of completely solvable and of rigid type if and
only if $\mathfrak g$ is nilpotent, that is, $M$ is a nilmanifold.
\item[ (2)]
A hyperelliptic surface can be characterized as
a 4-dimensional solvmanifold of rigid type with K\"ahler structure [\ref{H2}].
\end{list}
\bigskip

%%%%%%%%%%Example 3
\noindent{\bfseries Example 3}.
Let $G =: \C^l \rtimes \R^{2k}$, where the action
$\phi: \R^{2k} \rightarrow {\rm Aut}(\C^l)$ is defined by
$$\phi(\bar{t}_i) ((z_1, z_2, \ldots, z_l)) = 
(e^{\sqrt{-1}\,\eta^i_1\, t_i} z_1, e^{\sqrt{-1}\,\eta^i_2\, t_i} z_2, ...,
e^{\sqrt{-1}\, \eta^i_l\, t_i} z_l),$$
where $\bar{t}_i= t_i e_i$ ($e_i$: the $i$-th unit vector in $\R^{2k})$, and
$e^{\sqrt{-1}\,\eta^i_j}$ is the $s_i$-th root of unity,
$i = 1, \ldots,2k, j = 1, \ldots,l$. 
\medskip

If an abelian lattice $\Z^{2l}$ of $\C^l$ is preserved by
the action $\phi$ on $\Z^{2k}$, then $M=\Gamma \backslash G$ defines a solvmanifold
of rigid type, where $\Gamma=\Z^{2l} \rtimes \Z^{2k}$ is a lattice of $G$.
\bigskip

%%%%%%%%%%%%%%
The Lie algebra $\mathfrak g$ of $G$ is the following:
$$ \mathfrak g = \{X_1, X_2, \ldots , X_{2l}, X_{2l+1}, \ldots , X_{2l+2k}\}_\R,$$
where the bracket multiplications are defined by 
$$[X_{2l+2i}, X_{2j-1}] = -X_{2j}, [X_{2l+2i}, X_{2j}] = X_{2j-1}$$
for $i = 1, \ldots,k, j = 1, \ldots, l$,
and all other brackets vanish.
\medskip

 The canonical left-invariant complex structure is defined by
$$JX_{2j-1}=X_{2j}, JX_{2j}=-X_{2j-1},$$
$$JX_{2l+2i-1}=X_{2l+2i}, JX_{2l+2i}=-X_{2l+2i-1}$$ 
for $i=1, \ldots, k, j=1, \ldots, l$.
\medskip

%%%%%%%%%%%Example 4
\noindent{\bfseries Example 4}.
Let $G =: \C^l \rtimes \R^{2k}$, where the action
$\phi: \R^{2k} \rightarrow {\rm Aut}(\C^l)$ is defined by
$$\phi(\bar{t}_i) ((z_1, z_2, \ldots, z_l)) = 
(e^{2 \pi \sqrt{-1}\, t_i} z_1, e^{2 \pi \sqrt{-1}\, t_i} z_2, ...,
e^{2 \pi \sqrt{-1}\, t_i} z_l),$$
where $\bar{t}_i = t_i e_i$ ($e_i$: the $i$-the unit vector in $\R^{2k})$, 
$i = 1, \ldots, 2k$. 
\smallskip

Then, $\Z^{2n} \backslash G$
is a solvmanifold diffeomorphic to a torus $T^{2n}$ ($n= k+l$).
\bigskip

A compact solvmanifold $M$ in Example 3 is a finite quotient of a complex torus and
has a structure of a complex torus bundle over a complex torus, admiting
a canonical K\"ahler structure. We could have shown the converse
that if a compact solvmanifold admits a K\"ahler structure, then it must be of this type:
\medskip

%%%%%%%%%%%Theorem 5
\noindent{\bfseries Theorem 5} ([\ref{H2}, \ref{H3}]).
%\fbox{\begin{minipage}{18cm}
{\em A compact solvmanifold admits a
K\"ahler structure if and only if it is a finite quotient of a complex torus which
has a structure of a complex torus bundle over a complex torus.}
%\end{minipage}
%}
\medskip

We note that
\begin{list}{}{\topsep=5pt \leftmargin=12pt \itemindent=5pt \parsep=0pt \itemsep=5pt}
\item[ (1)] 
Since K\"ahlerian solvmanifolds (as defined in Example 3) are of rigid type, it follows
that a compact solvmanifold of completely solvable type has a K\"ahler
structure if and only if it is a complex torus. This is the so-called Benson-Gordon conjecture
([\ref{BG2}]).
\item[ (2)]  
We know [\ref{BG1}, \ref{H1}] that a compact nilmanifold
admits a K\"ahler structure if and only if it is a complex torus; and this result holds also for
bimeromorphic K\"ahler structures [\ref{H1}]. We see that
Theorem 5  also holds for bimeromorphic K\"ahler structures, since the proof is based on
this result and a result of Arapura and Nori [\ref{AN}] that a polycyclic K\"ahler group
must be almost nilpotent.
\item[ (3)] 
As noted in the paper [\ref{H3}], the Benson-Gordon conjecture
(stated in (1))
can be proved directly from the above results on K\"ahlerian nilmanifolds and
polycyclic K\"ahler groups, together with
a result of Auslander [\ref{AU}] that for a compact solvmanifold
$\Gamma \backslash G$, the Lie algebra $\mathfrak g$ of $G$ is of rigid type
if and only if $\Gamma$
is almost nilpotent (where $G$ is a simply connected solvable Lie group with discrete
subgroup $\Gamma$):
If $M = \Gamma \backslash G$ admits a K\"ahler structure and $\mathfrak g$ is of completely
solvable type, then $\Gamma$ is almost nilpotent. Hence  $\mathfrak g$ is both of rigid
type and of completely solvable type; and thus $\mathfrak g$ is nilpotent. Therefore,
$M$ is a compact K\"ahlerian nilmanifold, that is,  a complex torus.
There is a recent paper (by Baues and Cort\'es [\ref{BC}]) discussing
the Benson-Gordon conjecture and other relevant topics from more topological
point of view.

\end{list}
\bigskip

Concerning K\"ahler structures on a compact homogeneous manifold of compact
semi-simple Lie group, we have a fundamental theorem:
\medskip

%%%%%%%Theorem 6
\noindent{\bfseries Theorem 6} (Borel [\ref{B}], Goto [\ref{G}]).
%\fbox{\begin{minipage}{18cm}
{\em Let $G$ be a compact real semi-simple Lie group and
$D$ is a closed subgroup  which is the centralizer of a toral subgroup of $G$.
Then, $M = D \backslash G$ (of even-dimension) admits a homogeneous K\"ahler 
structure, which is a simply connected and projective
algebraic manifold. Conversely, if a compact homogeneous manifold of compact
semi-simple Lie group admits a K\"ahler structure, it  must be of the above form,
admitting a homogeneous K\"ahler structure.
}

%\end{minipage}
%}
\medskip

We note that

\begin{list}{}{\topsep=5pt \leftmargin=15pt \itemindent=5pt \parsep=0pt \itemsep=5pt}
\item[ (1)] 
$M = D \backslash G$ has a homogeneous complex structure $P \backslash G_\C$,
where $G_\C$ is the complexification of $G$, and $P$ is a parabolic subgroup of
$G_\C$ which contains a Borel subgroup $B$ of $G_\C$.
\item[ (2)] 
It is known (Samelson [\ref{S}], Wang [\ref{W}]) that any even-dimensional
compact semi-simple Lie group admits a homogeneous complex structure but
no K\"ahler structures.
\item[ (3)] It is known (Burstall et al. [\ref{BMGR}]) that  if a compact inner Riemannian
symmetric manifold admits a Hermitian structure (which is compatible with the given
metric), then it is Hermitian
symmetric. In particular, $S^6$ (considered as a compact inner Riemannian symmetric
manifold) admits no complex structures compatible with the given metric.

\end{list}
\bigskip
%\vspace{0.4in}

%%%%%%%%%%%%%%%%%%%%%%%%%%%%%%%%%%%%%%%
\noindent{\large \em \bfseries 3. The classification of 3-dimensional
compact complex solvmanifolds}
\bigskip

\noindent Let $M$ be a 3-dimensional compact complex solvmanifold.  Then,
$M$ can be written as $\Gamma \backslash G$, where $\Gamma$ is a lattice of
a simply connected unimodular complex solvable Lie group $G$ (cf. [\ref{BO}]).
\medskip

The Lie algebra $\g$ of $G$ is unimodular (i.e. the trace of ${\rm ad}\, (X) =0$
for every $X$ of $\g$), which is one of the following types:
\medskip

\begin{list}{}{\topsep=5pt \leftmargin=15pt \itemindent=5pt \parsep=0pt \itemsep=0pt}
\item[ (1)] Abelian Type: \par
$[X, Y] = [Y, Z] = [X, Z] = 0$.
\item[ (2)] Nilpotent Type: \par
$[X, Y] = Z,\; [X, Z] = [Y, Z] = 0$.
\item[ (3)] Non-Nilpotent Type: \par
$[X, Y] = -Y,\; [X, Z] = Z,\; [Y, Z] = 0$.
\end{list}
\medskip

\begin{list}{}{\topsep=5pt \leftmargin=15pt \itemindent=5pt \parsep=0pt \itemsep=5pt}

%1
\item[(1)] {\em \bfseries Abelian Type}:\; $G = \C^3$ \par
A lattice $\Gamma$ of $G$ is generated by a basis of $\C^3$ as a vector space
over $\R$.
\bigskip

%%%%%%%%%%%%
%2
\item[ (2)] {\em \bfseries Nilpotent Type}:\; $G = \C^2 \rtimes \C$ with
the action $\phi$ defined by
$$\phi(x)(y, z) = (y, z + xy),
$$
or in the matrix form,
$$G = \Bigg\{ \left(
\begin{array}[c]{ccc}
1 & x & z\\
0 & 1 & y\\
0 & 0 & 1
\end{array}
\right) \rule[-7mm]{0.25mm}{16mm}\;  x, y, z \in \C \Bigg\}.
$$
A lattice $\Gamma$ of $G$ can be written as 
$$\Gamma = \Delta \rtimes \Lambda,
$$
where $\Delta$ is a lattice of $\C^2$ and $\Lambda$ is a lattice of $\C$.
\medskip

Since an automorphism $f \in {\rm Aut}(\C)$ defined by
$f(x) = \alpha x, \alpha \not= 0$ can be extended to
an automorphism $F \in {\rm Aut}(G)$ defined by $F(x, y, z) = (\alpha x, \alpha^{-1} y, z)$,
we can assume that $\Lambda$ is generated by $1$ and
$\lambda\; (\lambda \notin \R)$ over $\Z$.
\medskip

 Since $\Delta$ is preserved by $\phi(1)$ and $\phi(\lambda)$, we see that
$\Delta$ is generated by
$(\alpha_1, \beta_1), (\alpha_2, \beta_2), (0, \alpha_1), (0, \alpha_2)$ over $\Z$,
where $\beta_1$ and $\beta_2$ are arbitrary complex numbers, and $\alpha_1$ and
$\alpha_2$ are linearly independent over $\R$ such
that $(\alpha_1, \alpha_2)$ is an eigenvector of some 
$A \in {\rm GL}(2, \Z)$ with the eigenvalue $\lambda$.
\medskip

Conversely, for any $A \in {\rm GL}(2, \Z)$ with non-real eigenvalue $\lambda$,
we can define a lattice $\Gamma$ of $G$. 
\medskip

\item[ (3)] {\em \bfseries Non-Nilpotent Type}:\; $G = \C^2 \rtimes \C$ with
the action $\phi$ defined by
$$\phi(x)(y, z) = (e^{x} y, e^{-x} z),
$$
or in the matrix form,
$$G = \Bigg\{ \left(
\begin{array}[c]{cccc}
e^{x} & 0 & 0 & y\\
0 & e^{-x} & 0 & z\\
0 & 0 & 1 & x\\
0 & 0 & 0 & 1
\end{array}
\right) \rule[-7mm]{0.25mm}{16mm}\;  x, y, z \in \C \Bigg\}.
$$
\medskip

%%%%%%%%%%%%%%
A lattice $\Gamma$ of $G$ can be written as $\Gamma = \Delta \rtimes \Lambda$,
where $\Delta$ is a lattice of $\C^2$, and $\Lambda$ is a lattice of $\C$
generated by $\lambda$ and $\mu$ over $\Z$. 
\medskip

Since $\Delta$ is preserved by $\phi(\lambda)$ and $\phi(\mu)$,
we see that $\Delta$ is generated by $(\alpha_i, \beta_i), i = 1, 2, 3, 4$ over $\Z$
such that 
$$\textstyle \gamma^{-1} \alpha_i = \sum_{j=1}^{4} a_{i j} \alpha_j,\; 
\gamma \beta_i = \sum_{j=1}^{4} a_{i j} \beta_j,
$$
\vspace{-10pt}
$$\textstyle \delta^{-1} \alpha_i = \sum_{j=1}^{4} b_{i j} \alpha_j,\; 
\delta \beta_i = \sum_{j=1}^{4} a_{i j} \beta_j,
$$
where $\gamma = e^{\lambda}, \delta = e^{\mu}$, and
$A = (a_{i j}), B = (b_{i j}) \in {\rm SL}_4 (\Z)$ are semi-simple and
mutually commutative.
\medskip

In other words, we have simultaneous eigenvectors 
$\alpha = (\alpha_1, \alpha_2, \alpha_3, \alpha_4),
\beta = (\beta_1, \beta_2, \beta_3, \beta_4) \in \C^4$ 
of $A$ and $B$ with eigenvalues $\gamma^{-1}, \gamma$ and $\delta^{-1}, \delta$
respectively.
\medskip

Conversely, for any mutually commutative, semi-simple matrices
$A, B \in {\rm SL}(4, \Z)$ with eigenvalues $\gamma^{-1}, \gamma$ and
$\delta^{-1}, \delta$ respectively, take simultaneous eigenvectors $\alpha, \beta \in \C^4$
of $A$ and $B$.
Then, $(\alpha_i, \beta_i), i = 1, 2, 3, 4$ are linearly independent over $\R$, defining
a lattice of $\Delta$ preserved by $\phi(\lambda)$ and $\phi(\mu)$ 
($\lambda = \log \gamma, \mu = \log \delta$).
\smallskip

Since $\lambda$ and $\mu$ are linearly independent over $\R$,
we have either $|\gamma| \not=1$ or $|\delta| \not= 1$. And if, for instance,
$|\gamma| \not=1$ and $\gamma \notin \R$,
then $A$ has four distinct eigenvalues
$\gamma^{-1}, \gamma, \overline{\gamma}^{-1}, \overline{\gamma}$.
\smallskip

For the case where both $A$ and $B$ have real eigenvalues
$\gamma^{-1}, \gamma$ and $\delta^{-1}, \delta$ respectively,
take simultaneous non-real eigenvectors $\alpha, \beta \in \C^4$ for them;
then we see that
$(\alpha_i, \beta_i), i = 1, 2, 3, 4$ are linearly independent over $\R$, defining
a lattice $\Delta$ of $\C^2$ preserved by $\phi(\lambda)$ and $\phi(\mu)$.
\bigskip
\end{list}

%%%%%%%%%%Example 5
\noindent{\bfseries Example 5}.
The Iwasawa manifold is obtained by
putting $\lambda = \sqrt{-1}, \alpha_1 = \alpha_2 = 0, \beta_1 = 1, \beta_2 = \sqrt{-1}$. 
\bigskip

%%%%%%%%%Example 6
\noindent{\bfseries Example 6}.
Take $A \in {\rm SL} (4, \Z)$ with four non-real eigenvalues
$\gamma, \gamma^{-1}, \overline{\gamma}, \overline{\gamma}^{-1}$;
for instance,
$$A = \left(
\begin{array}[c]{cccc}
0 & 1 & 0 & 0\\
0 & 0 & 1 & 0\\
0 & 0 & 0 & 1\\
-1 & 1 & -3 & 1
\end{array}
\right),
$$
with the characteristic polynomial given by
$${\rm det} (t I - A) = t^4 - t^3 + 3 t^2 - t + 1.$$
For the lattice $\Lambda$ of $\C$ generated by
$\lambda\;(\lambda = \log \gamma)$ and
$\mu = k \pi \sqrt{-1}\;(k \in~\Z)$,
and the lattice $\Delta$ of $\C^2$ generated by
$(\alpha_i, \beta_i), i = 1, 2, 3, 4$,
we can define a lattice $\Gamma = \Delta \rtimes \Lambda$ of $G$,
where $(\alpha_1, \alpha_2, \alpha_3, \alpha_4),
(\beta_1, \beta_2, \beta_3, \beta_4) \in \C^4$ are eigenvectors
of $A$ with eigenvalue $\gamma, \gamma^{-1}$.
\bigskip

%\break
%%%%%%%%%%Example 7
\noindent{\bfseries Example 7} (Nakamura [\ref{N}]).
Take $A \in {\rm SL}_2 (\Z)$ with two real
eigenvalues $\gamma^{-1}, \gamma$, $\gamma \not= \pm 1$, and
their real eigenvectors 
$(a_1, a_2), (b_1, b_2) \in \R^2$. 
Then, for any $\epsilon \notin \R$ (e.g. $\epsilon = \sqrt{-1}$),
$(a_1, a_2, a_1 \epsilon, a_2 \epsilon)$
and $(b_1, b_2, b_1 \epsilon, b_2 \epsilon)$ are
non-real eigenvectors for 
$A \oplus A \in {\rm SL}_4 (\Z)$ with eigenvalues $\gamma^{-1}, \gamma$.

For the lattice $\Lambda$ of $\C$ generated by
$\lambda\;(\lambda = \log \gamma)$ and
$\mu = k \pi \sqrt{-1}\;(k \in \Z)$,
and the lattice $\Delta$ of $\C^2$ generated by $(a_1, b_1),
(a_2, b_2), (a_1 \epsilon, b_1 \epsilon), (a_2 \epsilon, b_2 \epsilon)$,
we can define a lattice $\Gamma = \Delta \rtimes \Lambda$ of $G$.
%\end{list}

\bigskip

%%%%%%%%%%Theorem 7
\noindent{\bfseries Theorem 7} (Winkelmann [\ref{Wi}]).
%\fbox{\begin{minipage}{18cm}
{\em Let $G$ be a simply connected complex solvable linear algebraic group
with lattice $\Gamma$. Then, we have
$${\rm dim}\, H^1(\Gamma \backslash G, {\mathcal O}) =
{\rm dim}\, H^1(\g, \C) + {\rm dim}\, W,
$$
where ${\mathcal O}$ denotes the structure sheaf of $M$, $\n$ the nilradical of $\g$,
and $W$ the maximal linear subspace of $[\g, \g]/[\n, \n]$ for which ${\rm Ad}(\xi)$ on
$W$ is a real semi-simple linear endomorphism
for any $\xi \in \Gamma$. }
%\end{minipage}
%}
\bigskip 

%\underline{\bfseries Remark}\;
We have ${\rm dim}\, H^1(\g, \C) = {\rm dim}\,\g - {\rm dim}\,[\g, \g]$, and
${\rm Ad}(\xi)|W$ is diagonalizable over $\R$.
\medskip

Applying the Winkelmann's formula above and our classification of three-dimensional
compact complex solvmanifolds, 
we can determine $h^1(M) = {\rm dim}\, H^1(M, {\mathcal O})$ completely:

\begin{list}{}{\topsep=5pt \leftmargin=15pt \itemindent=5pt \parsep=0pt \itemsep=5pt}
\item[ (1)\,] Abelian Type:  ${\rm dim}\, W = 0$, $h^1 = 3$;
\item[ (2)\,] Nilpotent Type: ${\rm dim}\, W = 0$, $h^1 = 2$;
\item[ (3a)] Non-Nilpotent Type with $\gamma, \delta \in \R$:
${\rm dim}\, W = 2$, $h^1 = 3$;
\item[ (3b)] Non-Nilpotent Type with either $\gamma\; {\rm or}\; \delta \notin \R$:
${\rm dim}\, W = 0$, $h^1 = 1$.
\end{list}
\medskip 

\noindent{\bfseries Example 8}.
We see that Example 4 is of type (3a), and Example 5 is of type (3b).
\bigskip

%\vspace{0.5in}
\bigskip
%%%%%%%%%%%%%%%%%%%%%%%%%%%%%%%
\noindent{\large \em \bfseries 4. Pseudo-K\"ahler structures on
a $3$-dimensional compact complex solvmanifold}
\bigskip

We can see from Theorem 3 that a compact solvmanifold admits a homogeneous
pseudo-K\"ahler if and only if it is a complex torus. Yamada gave the first example of 
homogeneous complex pseudo-K\"ahler non-toral solvmanifold; and showed
the following fundamental result:
\bigskip

\noindent 
%%%%%%%%%%Theorem 8
\noindent{\bfseries Theorem 8} (Yamada [\ref{Y1}, \ref{Y2}]).
%\fbox{\begin{minipage}{18cm}
Let $M$ be an $n$-dimensional compact complex solvmanifold
which admits a pseudo-K\"ahler structure. Then,
we have $h^1(M)  \ge n$; and $M$ has a structure of a complex torus bundle over a
complex torus.

%\end{minipage}
%}
\medskip

%\underline{\bfseries Remark}\;
We remark that Winkelmann's formula implies that if  we have $h^1  \ge n$ then $[\n, \n] = 0$;
and thus the Mostow fibration gives a structure of a complex torus bundle over a
complex torus.
\bigskip

%%%%%%%Theorem 9
\noindent{\bfseries Theorem 9} ([\ref{H4}]).
%\fbox{\begin{minipage}{18cm}
{\em A 3-dimensional compact complex solvmanifold $M$
admits a pseudo-K\"ahler
structure if and only if it is of abelian type, or of non-nilpotent type with
$\gamma, \delta \in \R$.}
%\end{minipage}
%}
\medskip

\noindent{\bfseries Proof} (Sketch).
\medskip

If $M$ is of type (2) or (3b), then $M$ admits no pseudo-K\"ahler structures.
Therefore, it suffices to show that $M$ of type (3a) admits a pseudo-K\"ahler
structure.
\medskip

%%%%%%%%%%%%%%%%%%
We have $\gamma, \delta \in \R$ if and only if
$\Lambda$ is generated by $\lambda = a + k \pi \sqrt{-1}, \mu = b + l \pi \sqrt{-1}$,
where $a, b \in \R\; {\rm and}\; k, l \in \Z$.
\medskip

We can construct a
pseudo-K\"ahler structure $\omega$ on $\Gamma \backslash G$ in the following:
$$ \omega = \sqrt{-1} dx \wedge d \overline{x} + dy \wedge d \overline{z}
+ d \overline{y} \wedge d z,
$$
or using Maure-Cartan forms,
$\omega_1, \omega_2, \omega_3$,on $G$,
$$ \omega = \sqrt{-1} \omega_1 \wedge \overline{\omega_1} +
e^{- 2\, {\rm Im} (x) \sqrt{-1}}\, \omega_2 \wedge \overline{\omega_3} +
e^{2\, {\rm Im} (x) \sqrt{-1}}\, \overline{\omega_2} \wedge \omega_3,
$$
where $\omega_1 = dx, \omega_2 = e^x\, dy, \omega_3 = e^{-x}\, dz$.
\hfill $\Box$
\bigskip

Concerning pseudo-K\"ahler structures on compact complex nilmanifolds, we have
\medskip

%%%%%%%Theorem 10
\noindent{\bfseries Theorem 10} (Kodaira [\ref{N}]).
%\fbox{\begin{minipage}{18cm}
{\em Let $M$ be an $n$-dimensional compact complex nilmanifold,
and denote by $r$ the number of linearly independent closed holomorphic
$1$-forms on $M$. Then, we have $h^1(M) = r$, and $r=n$ holds if and only if
$M$ is a complex torus.}
%\end{minipage}
%}
\medskip

In particular, applying Theorem 8, we see that a non-toral compact complex nilmanifold admits
no pseudo-K\"ahler structures.
\medskip

We have the following result on
holomorphic principal fiber bundles over a complex torus with fiber a complex torus: 
\medskip

%%%%%%%Theorem 11
\noindent{\bfseries Theorem 11} (Murakami [\ref{Mu}]).
%\fbox{\begin{minipage}{18cm}
{\em A holomorphic principal fiber bundle over a complex torus
with fiber a complex torus  is a compact $2$-step nilmanifold
with a left-invariant complex structure: and it has a holomorphic connection
if and only if it is a compact complex nilmanifold.
}
%\end{minipage}
%}
\medskip

We see in particular that a holomorphic principal bundle over a complex torus
with fiber a complex torus admits no pseudo-K\"ahler structures.
\bigskip

%\vspace{0.5in}
%%%%%%%%%%%%%%%%%%%%%%%%%%%%%%%%%%%%%%%%
\noindent{\large \em \bfseries 5. Small deformations and non-left invariant
complex structures on a compact complex solvmanifold}
\bigskip

\noindent Let $G$ be a connected simply connected Lie group of dimension $2m$, and
$\g$ the Lie algebra of $G$.
\bigskip

\noindent{\bfseries Lemma 1}.
%\fbox{\begin{minipage}{18cm}
{\em  An almost complex structure $J$ on $\g$
is integrable
if and only if the subspace $\h_J$ of $\g_\C$ generated by $X + \sqrt{-1} JX  (X \in \g)$
is a complex subalgebra of $\g_\C$ such that $\g_\C = \h_J \oplus \overline{\h_J}$.}
%\end{minipage}
%}
\medskip

\noindent{\bfseries Lemma 2}.
%\fbox{\begin{minipage}{18cm}
 {\em Let $V$ be a real vector space of dimension $2m$.
Then, for a complex subspace $W$ of $V \otimes \C$ such that
$V \otimes \C = W \oplus \overline{W}$, there exists a unique $J_W \in {\rm GL}(V, \R),
{J_W}^2 = -I$ such that $W = \{X + \sqrt{-1} J_W X| X \in V \}_\C$.} 
%\end{minipage}
%}
\medskip

There exists one to one correspondence between complex structures $J$ on $\g$
and complex Lie subalgebras $\h$ such that $\g_\C = \h \oplus \overline{\h}$,
given by
$J \rightarrow \h_J$ and $\h \rightarrow J_\h$.
\bigskip

For a complex structure $J$, the complex Lie subgroup
$H_J$ of $G_\C$ corresponding to $\h_J$ is closed, simply connected,
and $H_J \backslash G_\C$ is biholomorphic to $\C^m$. The canonical inclusion
$\g \hookrightarrow \g_\C$ induces an inclusion $G \hookrightarrow G_\C$,
and $\Gamma = G \cap H_J$ is a discrete subgroup of $G$.
We have the following canonical map $g = i \circ \pi$:
$$G \stackrel{\pi}{\rightarrow} \Gamma \backslash G \stackrel{i}{\hookrightarrow}
H_J \backslash G_\C,
$$
where $\pi$ is a covering map, and $i$ is an inclusion. 
The left-invariant complex structure $J$ on $G$ is the one induced by $g$ from 
an open set $U = {\rm Im}\, g \subset \C^m$. For the details of the above argument
we refer to the paper [\ref{SN}].
\bigskip

%\vspace{0.4in}
%%%%%%%%%%%%%%%%%

Let $G$ be a 3-dimensional complex solvable Lie group of non-nilpotent
type, and $\g$ its Lie algebra.
Recall that $\g$ has a basis $X, Y, Z$ over $\C$ with bracket multiplication defined by
$$[X, Y] = -Y,\; [X, Z] = Z,\; [Y, Z] = 0.
$$

Let $\g_\R$ denote the real Lie algebra underlying $\g$, and
$\g_\C$ the complexification of $\g_\R$, that is,
$\g_\C = \g_\R \oplus \sqrt{-1} \g_\R$.
\medskip

Let $J_0$ be the original complex structure with
its associated complex subalgebra $\h_0$ of $\g_\C$
such that $\g_\C = \h_0 \oplus \overline{\h_0}$, and $H_0$
the complex subgroup of $G_\C$ corresponding to $\h_0$.
\bigskip

%%%%%%%%%Lemma 3
\noindent{\bfseries Lemma 3} ([\ref{H4}]).
%\fbox{\begin{minipage}{18cm}
{\em For any complex structure $J$ on $G$ with its
associated complex subalgebra 
$\h$ of $\g_\C$ such that $\g_\C = \h \oplus \overline{\h}$,
there exists a complex automorphism of
Lie algebras $\Phi:\g_\C \rightarrow \g_\C$ such that
$\Phi \circ \tau_0 = \tau \circ \Phi$ and $\Phi (\h_0) = \h$,
where $\tau_0$ and $\tau$ are the conjugations
with respect to $J_0$ and $J$ respectively.}
%\end{minipage}
%}
\smallskip

As a consequence we have
\medskip

%%%%%%%%%Theorem 12
\noindent{\bfseries Theorem 12} ([\ref{H4}]).
%\fbox{\begin{minipage}{18cm}
{\em Let $G$ be  a $3$-dimensional simply connected
complex solvable Lie group of non-nilpotent type. Then,
any left-invariant complex structure on $G$ is biholomorphic to $\C^3$.}
%\end{minipage}
%}
\medskip

\noindent{\bfseries Proof}.
The complex automorphism of Lie algebras
$\Phi$ induces a complex automorphism of Lie group
$\Psi: G_\C \rightarrow G_\C$ such that $q \circ \Psi= \tilde{\Psi} \circ q_0$,
which send $H_0$ to $H$ biholomorphically;
$$\begin{array}{ccccc}
(G, J_0) & \stackrel{i}{\hookrightarrow} & G_\C & \stackrel{q_0}{\rightarrow} &
H_0 \backslash G_\C \\
 &  & \vcenter{\rlap{$\scriptstyle \Psi$}}\: \downarrow & &
\vcenter{\rlap{$\scriptstyle \tilde{\Psi}$}}\: \downarrow \\
(G, J) & \stackrel{i}{\hookrightarrow} & G_\C & \stackrel{q}{\rightarrow} &
H \backslash G_\C 
\end{array}
$$
Here, we have $\Gamma = G \cap H_0 = \{0\}$, 
and $g_0 = q_0 \circ i$ is a biholomorphic map to $H_0 \backslash G_\C = \C^3$.
\hfill $\Box$
\bigskip

%%%%%%%%%%%%%%%%%
 We can also see that all left-invariant complex structures on
a $3$-dimensional simply connected complex solvable Lie group
are biholomorphic to $\C^3$.
\bigskip

%\vspace{0.4in}

Nakamura constructed small deformations of  $3$-dimensional
compact complex solvmanifolds;
and showed in particular that there exists a continuous family of complex
structures on those of type (3b) whose universal coverings are not Stein
(as noted in the paper, this construction is actually due to Kodaira). 
\bigskip

%\break
%%%%%%%%Theorem 13
\noindent{\bfseries Theorem 13} ([\ref{H4}]).
%\fbox{\begin{minipage}{18cm}
{\em There exists a continuous family of non-left-invariant complex structures on
the 6-dimensional compact solvmanifold $M$.}
%\end{minipage}
%}
\medskip

%\underline{\bfseries Remark}\;

We note that small deformations of a $3$-dimensional compact complex nilmanifold
(Iwasawa manifold) are all left-invariant (due to Salamon [\ref{S}]).
We conjecture that this also holds for
higher dimension. Recently, there appears a preprint (by McLaughlin et al. [\ref{MPPS}])
which proves the conjecture for a more general class of left-invariant complex 
structures on compact nilmanifolds.
\bigskip
%\break

\noindent {\bfseries Acknowledgements}\\
The author would like to thank the referee for his carefully reading the manuscript
with many useful comments and remarks.
\medskip

\begin{center}
{\bfseries References}
\end{center}
\renewcommand{\labelenumi}{[\theenumi]}
\begin{enumerate}
\itemsep=0pt
\item \label{AN} D. Arapura and M. Nori: {\em Solvable fundamental groups of
algebraic varieties and K\"ahler manifolds}, Compositio Math., {\bf 116} (1999), 173-188. 
\item \label{AU} L. Auslander: {\em An exposition of the structure of solvmanifolds I II}\,,
Bull. Amer. Math. Soc., {\bf 79} (1973),
No. 2, 227-261, 262-285.
\item \label{BC} O. Baues and V. Cort\'es: {\em Aspherical K\"ahler manifolds with solvable
fundamental group},  Geom. Dedicata, {\bf 122} (2006), 215-229.
\item \label{B} A. Borel: {\em K\"ahlerian coset spaces of semi-simple Lie groups},
Proc. Nat. Acad. Sci., {\bf 40} (1954), 1147-1151.
\item \label{BG1} C. Benson and C. S. Gordon: {\em K\"ahler and symplectic
structures on nilmanifolds}, Topology, {\bf 27} (1988), 755-782.
\item \label{BG2} C. Benson and C. S. Gordon: {\em K\"ahler structures on
compact solvmanifolds}, Proc. Amer. Math. Soc., {\bf 108} (1990), 971-980. 
\item \label{BO} W. Barth and M. Otte: {\em \"Uber fast-umiforme Untergruppen
komplexer Liegruppen und aufl\"osbare komplexe Mannigfaltigkeiten},
Comment. Math. Helv., {\bf 44} (1969), 269-281.
\item \label{BMGR} F. Burstall, O.Muskarov, G. Grantcharov and J. Rawnsley:
{\em Hermitian structures on Hermitian symmetric spaces},
J.  Geom. and Physics, {\bf 10} (1993) 245-249 .
\item \label{BR} A. Borel and R. Remmert: {\em \"Uber kompakte homogene K\"ahlersche
Mannigfaltigkeiten}, Math. Ann., {\bf 145} (1962), 429-439. 
\item \label{DG} J. Dorfmaister and D. Guan: {\em Classification of compact
pseudo-K\"ahler manifolds}, Comment Math. Helv., {\bf 67} (1992), 499-513.
\item \label{G} M. Goto: {\em On algebraic homogeneous spaces}, Am. J. Math.,
{\bf 76} (1954), 811-818.
\item \label{Gu} D. Guan: {\em Modification and cohomology groups of compact
solvmanifolds}, Preprint.
\item \label{H1} K. Hasegawa: {\em Minimal models of nilmanifolds},
Proc. Amer. Math. Soc., {\bf 106} (1989), 65-71.
\item \label{H2} K. Hasegawa: {\em Complex and K\"ahler structures on
compact solvmanifolds}, Proceedings of the conference on symplectic topology,
Stare Jablonki (2004), J. Symplectic Geom., {\bf 3} (2005), 749-767.
\item \label{H3} K. Hasegawa: {\em A note on compact solvmanifolds with
K\"ahler structures}, Osaka J. Math., {\bf 43} (2006), 131-135.
\item \label{H4} K. Hasegawa: {\em Small deformations and non-left-invariant
complex structures on a compact solvmanifold}, Preprint, math.CV/0703756.
\item \label{M} Y. Matsushima: {\em Sur les espaces homog\`enes K\"ahl\'eriens
d'un groupe de Lie r\'eductif}, Nagoya Math. J. {\bf 11} (1957), 53--60.
\item \label{MPPS} C. McLaughlin, H. Pedersen, Y. S. Poon and S. Salamon:
Deformation of 2-Step Nilmanifolds with abelian complex structures, Preprint.
\item \label{Mu} S. Murakami: {\em Sur Certains Espaces Fibr\'es Principaux
Holomorphes donts le Groupe est Ab\'elien Connexe}, Osaka J. Math.
{\bf 13} (1961), 143-167.
\item \label{N} I. Nakamura: {\em Complex parallelizable manifolds and their
small deformations}, J. Differential Geom., {\bf 10} (1975), 85-112.
\item \label{OR} K. Oeljeklaus and W. Richthofer: {\em Homogeneous complex
surfaces}, Math. Ann., {\bf 268} (1984), 273-292.
\item \label{S} S. M. Salamon: {\em Complex structures on nilpotent Lie algebras},
J. Pure and Appl. Algebra, {\bf 157}, 311-333.
\item \label{SA} H. Samelson: {\em A class of complex analytic manifolds},
Portugaliae Math., {\bf 12} (1953), 129-132.
\item \label{SN} D. Snow: {\em Invariant complex structures on reductive Lie groups},
J. Reine Angew. Math., {\bf 371} (1986), 191-215.
\item \label{T} J. Tits: {\em Espaces homog\`enes complexes compacts},
Comm. Math. Helv., {\bf 37} (1962), 111-120.
\item \label{W} H. C. Wang: {\em Closed manifolds with homogeneous complex structure},
Am. J. Math., {\bf 76} (1954), 1-32. 
\item \label{Y1} T. Yamada: {\em A pseudo-K\"ahler structure on a non-toral compact complex
parallelizable solvmanifold}, Geom. Dedicata, {\bf 112} (2005), 115--122.
\item \label{Y2} T. Yamada: {\em A structure theorem of compact complex parallelizable
pseudo-K\"ahler solvmanifolds}, Osaka J. Math., {\bf 43} (2006), 923-933.
\item \label{W} C.T.C. Wall: {\em Geometric structures on compact complex analytic surfaces},
Topology {\bf 25} (1986), 119--153.
\item \label{Wi} J. Winkelmann: Complex analytic geomtry of complex
parallelizable manifolds, M\'em. Soc. Math. de France S\'er. 2, 72-73 (1998).

\end{enumerate}

\begin{flushleft}
Department of Mathematics\\
Faculty of Education and Human Sciences\\
Niigata University, Niigata\\
JAPAN\\
\vskip3pt
e-mail: hasegawa@ed.niigata-u.ac.jp\\
\end{flushleft}

\end{document}